\documentclass[11pt,a4paper]{amsart}
\usepackage[T1]{fontenc}
\usepackage[utf8]{inputenc}
\usepackage{lmodern,microtype}
\usepackage{amsmath,amssymb,amsthm,mathtools}
\usepackage[margin=2.6cm]{geometry}
\usepackage{graphicx,tikz,enumitem,array}
\usetikzlibrary{arrows.meta,positioning,calc,decorations.pathreplacing}
\usepackage{xcolor,hyperref}
\hypersetup{colorlinks=true,linkcolor=black,citecolor=black,urlcolor=blue!45!black}
\newtheorem{theorem}{Theorem}[section]
\newtheorem{proposition}[theorem]{Proposition}
\newtheorem{lemma}[theorem]{Lemma}
\newtheorem{corollary}[theorem]{Corollary}
\theoremstyle{remark}
\newtheorem{remark}[theorem]{Remark}
\theoremstyle{definition}
\newtheorem{example}[theorem]{Example}
\newcommand{\dd}{\,\mathrm{d}}
\newcommand{\cS}{\mathcal S}
\newcommand{\wt}{\widetilde}
\newcommand{\tr}{\operatorname{tr}}

\title[Exit-pair rigidity and inverse reconstruction]{Exit-Pair Rigidity, Neutral Restart Laws, and Inverse Reconstruction in One Dimension}
\author{Tony Newton}
\address{Newton Astro Labs, United Kingdom}
\email{tony.newton79@gmail.com}
\author{Juan Antonio Vega Coso}
\address{IUFFyM, Universidad de Salamanca, Spain}
\email{jantovc@gmail.com}
\date{}
\begin{document}
\begin{abstract}
For an absorbed one-dimensional diffusion, we characterize two-point restart laws preserving the exit-side probability for every independent renewal clock. Spectral identities establish rigidity under spatial accumulation, while explicit nonsymmetric potentials admit neutral pairs. We distinguish midpoint Dirichlet isospectrality from Weyl reflection rigidity and prove a second-order high-energy asymptotic under local continuity, giving two necessary asymmetry conditions for neutral pairs. Two complete exit profiles at distinct reset rates identify diffusivity and drift; a finite birth--death analogue admits exact Jacobi criteria and reconstruction from four channel functions.
\end{abstract}
\maketitle
\pagestyle{plain}

\section*{Introduction}
Resetting can change which endpoint an absorbed process reaches, a phenomenon studied in the diffusion literature since Evans and Majumdar~\cite{EvansMajumdar}; see~\cite{EvansMajumdarSchehr} for a review, \cite{PalReuveni} for restart with general renewal clocks and~\cite{Belan} for the effect of restart on the probability of success. In a companion paper~\cite{NewtonVegaA} we showed, for nearest-neighbour walks, that the exit-side probability under renewal resetting depends on the walk only through the exit pair, the joint law of the absorption time $T$ and the exit side $J$, and that classical stochastic orders between the two conditional exit-time laws correspond to monotone responses to resetting. In particular, a site is reset-neutral for every renewal clock exactly when $T$ and $J$ are independent~\cite[Theorem~2.2]{NewtonVegaA}. Here the exit pair of a one-dimensional diffusion is taken as data, and we ask what it determines. The mean exit time alone is insufficient: the exit side and exit time jointly determine which reset schedules preserve a splitting probability. The same observables also carry spectral information. We investigate two consequences. First, neutrality of two-point restart laws can be rigid under spatial accumulation, even though individual neutral pairs exist for nonsymmetric potentials. Second, two full spatial splitting profiles at distinct Poisson reset rates determine both coefficients of a regular one-dimensional diffusion. Both questions are read from the same object, the exit channels $U_\pm(\lambda,y)$: as functions of the rate at a fixed starting point they carry spectral data, which govern neutrality, and as functions of the starting point at a fixed rate they are spatial profiles, which carry the coefficients.

The principal results are Theorem~\ref{thm:weights} (uniqueness of pair weights), Theorems~\ref{thm:interval}--\ref{thm:families} (symmetry and counterexamples), and Theorem~\ref{thm:two} (two-rate identification). The distinction between midpoint Dirichlet isospectrality and Weyl rigidity is developed in Section~\ref{sec:weyl}. Proposition~\ref{prop:single}, the continuous-time form of~\cite[Theorem~2.2]{NewtonVegaA}, is our starting point and is not claimed as new. The contributions are the neutral-pair rigidity and construction theorems in Sections~2--3, the midpoint Dirichlet/Weyl distinction in Section~4, the discrete Jacobi criterion and reconstruction (Theorem~\ref{thm:disc}, Proposition~\ref{prop:three}), and the two-rate identification (Theorem~\ref{thm:two}).

Relevant inverse first-passage work includes Bal--Chou~\cite{BalChou} and Belousov--Hassanali--Roldán~\cite{Belousov}, while reset-conditioned first passage and inverse reset problems have been studied by Pal--Prasad~\cite{PalPrasad}, Villarroel~\cite{Villarroel}, and Abundo~\cite{Abundo,Abundo2025}. Abundo studies recovery of the distribution of the initial position when the diffusion coefficients are known; here the coefficients themselves are unknown and the data are full spatial splitting profiles. Compared with Borg's two-spectra theorem and the Gel'fand--Levitan theory~\cite{Borg,GelfandLevitan}, which use spectral data at a fixed point, the identification here uses spatial profiles at a fixed spectral parameter. To our knowledge, recovery of both $D$ and $b$ from splitting profiles at two Poisson rates has not previously been established. Nickl and Seizilles~\cite{NicklSeizilles} study inference of diffusivity from killing locations, a different observation mechanism. The data used for identification are \emph{two full spatial splitting profiles}, not the single-start first-exit time law or a finite number of noisy measurements. The Schwarzian identity is classical~\cite{OvsienkoTabachnikov}; the identification theorem applies it to these particular exit-pair observables. Such reconstructions are of practical interest wherever splitting probabilities between two absorbing boundaries can be measured, as in the inference of position-dependent diffusivity and drift from first-passage statistics in inhomogeneous transport~\cite{Belousov}; resetting adds a controllable parameter, the reset rate, which supplies the second profile.

The results also admit an interpretation in terms of spectral duality structures~\cite{VegaSDS,VegaRigidity}. For a reflection-symmetric Liouville potential, the exit channels obey a mode-independent duality with a constant orbit product; the corresponding spectral relation follows by taking residues at the Dirichlet poles. The rigidity theorems distinguish this structural neutrality from isolated neutral pairs in nonsymmetric potentials. We develop this correspondence below without invoking the axioms of the finite-site theory.

\section{Exit-pair renewal identities and the intrinsic equation}\label{sec:framework}
\subsection{The observable and the renewal identity}
Let $T$ be the absorption time and $J\in\{0,L\}$ the endpoint reached. For a starting location $y$, define $\nu_y(\dd t)=\mathbb P_y(T\in\dd t,J=0)$ and $\eta_y(\dd t)=\mathbb P_y(T\in\dd t,J=L)$. An admissible independent renewal clock has independent, identically distributed inter-reset durations that are independent of the diffusion, with non-increasing right-continuous survival function $W$ and $W(0)=1$. After each reset the process restarts from a point drawn from a restart law $\mu$ on $(0,L)$, independently of everything else, and the initial point is drawn from $\mu$ as well; a fixed starting point $y$ is the case $\mu=\delta_y$, and general laws are used from Proposition~\ref{prop:laws} on. The clock is admissible when the probability of absorption before a reset is strictly positive. The resulting left-exit probability is
\begin{equation}\label{eq:renew}
 Q_y(W)=\frac{\int W(t)\nu_y(\dd t)}{\int W(t)[\nu_y+\eta_y](\dd t)}.
\end{equation}
Indeed, if $A$ and $B_0$ are the probabilities of left and right absorption before a reset in one cycle, then $Q=A+(1-A-B_0)Q$. For a Poisson clock of rate $\lambda$, set $U_-(\lambda,y)=\mathbb E_y[e^{-\lambda T};J=0]$, $U_+(\lambda,y)=\mathbb E_y[e^{-\lambda T};J=L]$ and $q_\lambda(y)=U_-/(U_-+U_+)$. Uniqueness of Laplace transforms yields the following equivalence, the diffusion form of~\cite[Theorem~2.2]{NewtonVegaA}.
\begin{proposition}[Universal individual neutrality]\label{prop:single}
For an interior starting point $y$ of a regular diffusion, the following are equivalent: $Q_y(W)$ is independent of every admissible independent renewal clock; $q_\lambda(y)$ is constant over $\lambda>0$; the random variables $T$ and $J$ are independent; and $\nu_y$ and $\eta_y$ are proportional finite measures.
\end{proposition}
\begin{proof}
Proportionality immediately gives clock independence in \eqref{eq:renew}. A Poisson clock of rate $\lambda>0$ is an admissible renewal clock with survival function $W(t)=e^{-\lambda t}$, for which \eqref{eq:renew} gives $Q_y(W)=q_\lambda(y)$. Hence clock independence implies constancy of $q_\lambda$. Conversely, constancy of $q_\lambda$ implies that the Laplace transform of $\nu_y-c(\nu_y+\eta_y)$ vanishes at every positive rate; uniqueness for finite signed measures gives proportionality, equivalent to independence of $T$ and $J$.
\end{proof}

\subsection{Physical and intrinsic coordinates}
A physical diffusion on $[0,\ell]$ has generator $\mathcal L_xf=D(x)f''+b(x)f'$ with $D\in C^2[0,\ell]$ strictly positive and $b\in C^1[0,\ell]$. Put
\begin{equation}\label{eq:lamperti}
 y=s(x)=\int_0^x\frac{\dd u}{\sqrt{D(u)}},\quad L=s(\ell),\quad
 \beta(s(x))=\frac{b(x)-D'(x)/2}{\sqrt{D(x)}}.
\end{equation}
Then $\mathcal L_y=\partial_y^2+\beta(y)\partial_y$, and
$b(x)=\sqrt{D(x)}\,\beta(s(x))+D'(x)/2$. The gauge $u=e^{-B/2}\psi$, $B(y)=\int_0^y\beta(t)\dd t$, gives the Sturm equation $\psi''=(\lambda+V)\psi$ with $V=\beta'/2+\beta^2/4$.
\subsection{The Sturm boundary problem and transfer matrices}

We work in the intrinsic coordinate $y\in[0,L]$, with generator $\partial_y^2+\beta\partial_y$, $B(y)=\int_0^y\beta$, $V=\beta'/2+\beta^2/4$, $\lambda=k^2>0$, Sturm solutions $S_\pm$ with $S_-(0)=0$, $S_-'(0)=1$, $S_+(L)=0$, $S_+'(L)=-1$, characteristic function $\Delta=S_-(L)=S_+(0)$, and define the channels
\begin{equation}\label{eq:U}
U_-(\lambda,y)=e^{-B(y)/2}\frac{S_+(\lambda,y)}{\Delta(\lambda)},\qquad
U_+(\lambda,y)=e^{[B(L)-B(y)]/2}\frac{S_-(\lambda,y)}{\Delta(\lambda)}.
\end{equation}
\textbf{Regularity.} Throughout, $\beta\in W^{1,1}(0,L)$. Then $B\in C^1$, $V\in L^1(0,L)$, and solutions of $\psi''=(\lambda+V)\psi$ are taken in the Carathéodory sense ($\psi,\psi'$ absolutely continuous); the initial value problem has a unique solution. By a \emph{regular diffusion} we mean one with generator $\partial_y^2+\beta\partial_y$ in the intrinsic coordinate, $\beta\in W^{1,1}(0,L)$, absorbed at $0$ and $L$; this includes the physical diffusions above and the drifts $\beta=2\psi_0'/\psi_0$ of Theorem~\ref{thm:real}. Such a diffusion is determined by its scale density $e^{-B}$ and speed density $e^{B}$, and $U_-$, $U_+$ in \eqref{eq:U} are the solutions of $u''+\beta u'=\lambda u$ with boundary values $(1,0)$ and $(0,1)$, that is, $U_-=\mathbb E_y[e^{-\lambda T};J=0]$ and $U_+=\mathbb E_y[e^{-\lambda T};J=L]$. Equalities between potentials, such as $V=\wt V$, periodicity or the conditions in Theorem~\ref{thm:families}, are understood almost everywhere. Write $M=\|V\|_{L^1}$, $R(y)=L-y$ and $\wt V=V\circ R$. We write $S_\pm(\lambda,y;V)$ when the dependence on the potential must be explicit.

\begin{lemma}[Positivity]\label{lem:pos}
For $\lambda>0$: $\Delta(\lambda)>0$, $S_->0$ on $(0,L]$ and $S_+>0$ on $[0,L)$. 
\end{lemma}
\begin{proof}
$f=e^{-B/2}S_-$ solves $f''+\beta f'=\lambda f$ with $f(0)=0$, $f'(0)=1$. If $y_1\in(0,L]$ were its first positive zero, then $f>0$ on $(0,y_1)$. The divergence-form identity $(e^B f')'=\lambda e^B f$ holds almost everywhere, so $e^B f'$ is increasing on $[0,y_1]$ and equals $1$ at zero. Thus $f'>0$ on this interval and $f(y_1)>f(0)=0$, a contradiction. The right-end argument is its reflection, and $\Delta=S_-(L)>0$.
\end{proof}

\begin{lemma}[Reversal]\label{lem:rev}
$S_+(\lambda,y;V)=S_-(\lambda,L-y;\wt V)$ for $y\in[0,L]$, and $\Delta(\lambda;V)=\Delta(\lambda;\wt V)$. In particular the Dirichlet characteristic function of an interval does not change when the potential is reversed. 
\end{lemma}
\begin{proof}
$\phi(s)=S_+(\lambda,L-s;V)$ solves $\phi''=(\lambda+\wt V)\phi$ with $\phi(0)=0$, $\phi'(0)=1$, so $\phi=S_-(\cdot;\wt V)$ by uniqueness. At $y=0$: $\Delta(V)=S_+(0;V)=S_-(L;\wt V)=\Delta(\wt V)$.
\end{proof}

\paragraph{Transfer matrices.} For a piece $X$ of the potential on an interval of length $\ell$, let $M_X=\begin{pmatrix}c&s\\c'&s'\end{pmatrix}(\ell)$, with $c,s$ the fundamental solutions ($c(0)=1$, $c'(0)=0$, $s(0)=0$, $s'(0)=1$); $\det M_X=1$ and $M_X$ maps $(\psi,\psi')$ at the left end to the right end. For a concatenation $X_1X_2\cdots X_n$ (left to right), $M_{X_1\cdots X_n}=M_{X_n}\cdots M_{X_1}$, and $S_-$ at the right end of the word is the $(1,2)$ entry, which we write $[\,\cdot\,]_{12}$. Let $E=\mathrm{diag}(1,-1)$ and $\wt X$ the reversed piece.

\begin{lemma}[Reversed piece]\label{lem:revmat}
$M_{\wt X}=E\,M_X^{-1}E$. Hence $[M_{\wt X}]_{12}=[M_X]_{12}$, and for any word, $M_{\wt X_n\cdots\wt X_1}=E\,(M_{X_1\cdots X_n})^{-1}E$. 
\end{lemma}
\begin{proof}
If $\psi$ solves the equation for $X$, then $\varphi(x)=\psi(\ell-x)$ solves it for $\wt X$, and $(\varphi,\varphi')(x)=E(\psi,\psi')(\ell-x)$. So $(\varphi,\varphi')(\ell)=E(\psi,\psi')(0)=EM_X^{-1}(\psi,\psi')(\ell)=EM_X^{-1}E(\varphi,\varphi')(0)$. For $M=\begin{pmatrix}a&b\\c&d\end{pmatrix}$ with $\det M=1$, $EM^{-1}E=\begin{pmatrix}d&b\\c&a\end{pmatrix}$. The statement for words follows from $(AB)^{-1}=B^{-1}A^{-1}$.
\end{proof}

\begin{lemma}[High-energy bound]\label{lem:volterra}
Let $w_-(y)=2k\,e^{-ky}S_-(\lambda,y)$. For all $y\in[0,L]$,
\[
\bigl|w_-(y)-(1-e^{-2ky})\bigr|\le\frac M{2k}\,e^{M/2k},\qquad |w_-(y)|\le e^{M/2k}.
\]
The same holds for $w_+(y)=2k\,e^{-ky}S_+(\lambda,L-y)$. Likewise $2e^{-k\ell}c(\ell)\to1$ and $2e^{-k\ell}s'(\ell)\to1$ for any piece, so $\tr M_X\sim e^{k\ell}$ as $k\to\infty$. 
\end{lemma}
\begin{proof}
Variation of constants gives $S_-(y)=\frac{\sinh ky}{k}+\int_0^y\frac{\sinh k(y-t)}{k}V(t)S_-(t)\dd t$, that is,
\begin{equation}\label{eq:volterra}
w_-(y)=1-e^{-2ky}+\frac1{2k}\int_0^y\bigl(1-e^{-2k(y-t)}\bigr)V(t)\,w_-(t)\dd t .
\end{equation}
Since $0\le1-e^{-2k(y-t)}\le1$, Gronwall gives $|w_-|\le e^{M/2k}$, and inserting this in the integral gives the first bound. For $w_+$ apply the same to $\wt V$ (Lemma~\ref{lem:rev}), which has the same $L^1$ norm. For $c$ and $s'$ the same Volterra argument applies to $c(x)=\cosh kx+\int_0^x\frac{\sinh k(x-t)}{k}Vc$ and to the derivative of the equation for $s$.
\end{proof}

\begin{lemma}[Order and normalization]\label{lem:order}
Let $W\in L^1(0,c)$ and $N=\|W\|_{L^1(0,c)}$. For each $x\in(0,c]$, $\lambda\mapsto s(\lambda,x;W)$ and $\lambda\mapsto s'(\lambda,x;W)$, with $s'=\partial_xs$, are entire functions of order at most $1/2$:
\[
|s(\lambda,x)|\le x\,e^{xN}e^{x|\lambda|^{1/2}},\qquad |s'(\lambda,x)|\le\bigl(1+Nxe^{xN}\bigr)e^{x|\lambda|^{1/2}}\qquad(\lambda\in\mathbb C).
\]
As $\lambda=k^2\to+\infty$,
\[
2k\,e^{-kx}s(\lambda,x)=1+o(1),\qquad 2e^{-kx}s'(\lambda,x)=1+o(1),\qquad \frac{s'(\lambda,x)}{s(\lambda,x)}=k+o(1).
\]
\end{lemma}

\begin{proof}
For $\lambda\in\mathbb C$ choose a square root $k$; the functions $\sinh(kx)/k$ and $\cosh(kx)$ are even in $k$, hence entire in $\lambda$. Since $\sinh(kx)/k=\int_0^x\cosh(kt)\dd t$, we have $|\sinh(kx)/k|\le x\,e^{|\mathrm{Re}\,k|x}$ and $|\cosh(kx)|\le e^{|\mathrm{Re}\,k|x}$. Variation of constants gives
\begin{equation}\label{eq:vc}
s(x)=\frac{\sinh kx}{k}+\int_0^x\frac{\sinh k(x-t)}{k}W(t)s(t)\dd t,\qquad
s'(x)=\cosh kx+\int_0^x\cosh k(x-t)\,W(t)s(t)\dd t .
\end{equation}
The Picard iterates of the first equation are entire in $\lambda$ and converge locally uniformly (the standard Volterra estimate), so $s(\cdot,x)$ is entire, and then so is $s'(\cdot,x)$ by the second equation. With $\varphi(x)=e^{-|\mathrm{Re}\,k|x}|s(x)|$ the first equation gives $\varphi(x)\le x+\int_0^x(x-t)|W(t)|\varphi(t)\dd t\le x+x\int_0^x|W|\varphi$, and Gronwall's inequality applied to $\varphi(x)/x$ gives $\varphi(x)\le x\,e^{xN}$. Inserting this in the second equation gives the bound for $s'$. Since $|\mathrm{Re}\,k|\le|\lambda|^{1/2}$, both functions have order at most $1/2$.

For $\lambda>0$ put $w(x)=2ke^{-kx}s(x)$ and $v(x)=2e^{-kx}s'(x)$. Multiplying \eqref{eq:vc} by $2ke^{-kx}$ and by $2e^{-kx}$,
\[
w(x)=1-e^{-2kx}+\frac1{2k}\int_0^x\bigl(1-e^{-2k(x-t)}\bigr)Ww\dd t,\qquad
v(x)=1+e^{-2kx}+\frac1{2k}\int_0^x\bigl(1+e^{-2k(x-t)}\bigr)Ww\dd t .
\]
As in Lemma~\ref{lem:volterra}, $|w|\le e^{N/2k}$, so $w(x)=1+O(1/k)$ and $v(x)=1+O(1/k)$. Moreover
\[
v(x)-w(x)=2e^{-2kx}+\frac1k\int_0^xe^{-2k(x-t)}W(t)w(t)\dd t=o(1/k),
\]
because $|e^{-2k(x-t)}W(t)w(t)|\le e^{N}|W(t)|$ for $k\ge1/2$ and the integrand tends to $0$ for almost every $t<x$ (dominated convergence). Hence $s'/s=k\,v/w=k\bigl(1+o(1/k)\bigr)=k+o(1)$.
\end{proof}

\begin{lemma}[Uniqueness from the Weyl function]\label{lem:weylunique}
Let $W_1,W_2\in L^1(0,c)$ be real-valued, write $s_i=s(\cdot,c;W_i)$ and $s_i'=s'(\cdot,c;W_i)$, and suppose that
\begin{equation}\label{eq:weyleq}
s_1'(\lambda)\,s_2(\lambda)=s_2'(\lambda)\,s_1(\lambda)
\end{equation}
for all $\lambda$ in a subset of $(0,\infty)$ with an accumulation point in $(0,\infty)$. Then $s_1\equiv s_2$, $s_1'\equiv s_2'$, and $W_1=W_2$ almost everywhere on $(0,c)$.
\end{lemma}

\begin{proof}
Both sides of \eqref{eq:weyleq} are entire in $\lambda$ (Lemma~\ref{lem:order}), so \eqref{eq:weyleq} holds on $\mathbb C$. For each $\lambda$, $s_i(\lambda)$ and $s_i'(\lambda)$ do not vanish together, otherwise the solution would vanish identically, against $s'(0)=1$. Let $\lambda_0$ be a zero of $s_1$ of multiplicity $r$. Then $s_1'(\lambda_0)\ne0$, so by \eqref{eq:weyleq} $\lambda_0$ is a zero of $s_2$ of multiplicity at least $r$; by symmetry $s_1$ and $s_2$ have the same zeros with the same multiplicities. Both have order at most $1/2<1$, so by Hadamard's factorization theorem each is $C_i\lambda^{m}\prod_n(1-\lambda/\lambda_n)$ over its non-zero zeros, with the same $m$ and the same product; hence $s_1=Cs_2$ for a constant $C\ne0$. By Lemma~\ref{lem:order}, $s_1/s_2\to1$ as $\lambda\to+\infty$, so $C=1$. Then \eqref{eq:weyleq} gives $s_1'=s_2'$ (divide by $s_1=s_2\not\equiv0$).

The zeros of $s_i$ are $\lambda=-\mu$ with $\mu$ in the Dirichlet spectrum of $-\dd^2/\dd x^2+W_i$ on $[0,c]$, and the zeros of $s_i'$ are $\lambda=-\mu$ with $\mu$ in its Dirichlet--Neumann spectrum (Dirichlet at $0$, Neumann at $c$). So $W_1$ and $W_2$ have the same two spectra, and Borg's two-spectra theorem~\cite{Borg}, in its real-integrable-potential form~\cite{LevitanGasymov} and~\cite[Theorem~4.4]{HrynivMykytyuk} (there the Neumann condition uses the quasi-derivative $y'-\sigma y$ with $\sigma'=W_i$; the choice $\sigma(x)=-\int_x^cW_i$ makes it $y'(c)$), gives $W_1=W_2$ almost everywhere.
\end{proof}

\subsection{Restart distributions and neutral laws}

\begin{proposition}[Neutral restart laws]\label{prop:laws}
Let $\mu$ be a probability measure on $(0,L)$, $U_\pm(\lambda,\mu)=\int U_\pm(\lambda,y)\,\mu(\dd y)$, and $\nu_\mu$, $\eta_\mu$ the corresponding averaged exit measures. The following are equivalent: the eventual left-exit probability is the same for every admissible independent renewal clock; $q_\lambda(\mu)=U_-(\lambda,\mu)/[U_-(\lambda,\mu)+U_+(\lambda,\mu)]$ is independent of $\lambda>0$; $\nu_\mu$ and $\eta_\mu$ are proportional. 
\end{proposition}
\begin{proof}
The renewal argument of Proposition~\ref{prop:single} applies because the restart point is redrawn independently at each renewal: with $A,B$ the one-cycle probabilities of left and right absorption, $Q=A/(A+B)$, now with $\nu_\mu,\eta_\mu$. The rest is the Laplace-uniqueness argument of Proposition~\ref{prop:single}.
\end{proof}

\begin{proposition}[Pair laws]\label{prop:pairform}
Let $0<y<L/2$, $w_1,w_2>0$ with $w_1+w_2=1$, $\mu=w_1\delta_y+w_2\delta_{L-y}$ and $C\in(0,1)$. Put
\[
a_1=w_1e^{-B(y)/2},\qquad a_2=w_2e^{-B(L-y)/2},\qquad c=\frac{C\,e^{B(L)/2}}{1-C}>0,\qquad \Phi=S_+-c\,S_- .
\]
Then $q_\lambda(\mu)=C$ for all $\lambda>0$ if and only if $a_1\Phi(\lambda,y)+a_2\Phi(\lambda,L-y)=0$ for all $\lambda>0$. 
\end{proposition}
\begin{proof}
$q_\lambda(\mu)=C$ is $(1-C)\,U_-(\lambda,\mu)=C\,U_+(\lambda,\mu)$. Insert \eqref{eq:U} and multiply by $\Delta/(1-C)$.
\end{proof}

\begin{remark}[Dictionary with spectral duality structures]\label{rem:sds}
The analogy with the finite spectral duality structure of~\cite{VegaSDS,VegaRigidity} has sites $y\in(0,L)$ and the involution $R(y)=L-y$, with unique fixed point $L/2$. The two channels are $U_\pm(\lambda,y)$; the spectral modes occur at the zeros of the Dirichlet characteristic function $\Delta(\lambda)$, understood through the meromorphic continuation of the channels. For a symmetric Liouville potential, Theorem~\ref{thm:real} establishes
\[
U_+(\lambda,y)=\kappa(y)U_-(\lambda,L-y),\quad
\kappa(y)>0,\quad\kappa(y)\kappa(L-y)=K=e^{B(L)}.
\]
The factor is independent of the rate and therefore also relates the residues mode by mode. The neutral law on a two-site orbit has the weights and common level given by Theorem~\ref{thm:weights}; constant drift $\beta\equiv\beta_0$ corresponds to the biased walk, with $K=e^{\beta_0L}$ in place of $(p/q)^a$. There are continuously many sites and infinitely many spectral modes, so this is a continuous analogue, not a finite structure in the axiomatic sense of~\cite{VegaSDS}. Theorems~\ref{thm:interval}--\ref{thm:endpoint} describe rigidity of neutrality under accumulation, while Theorem~\ref{thm:families} provides accidental nonsymmetric neutrality without a global reflection duality.
\end{remark}

\section{High-energy asymmetry and KdV densities}\label{sec:asym}
For $0<y<L$, write
\[
g(\lambda,y)=\frac{S_-(\lambda,y)}{S_+(\lambda,L-y)},\qquad
A(y)=\int_0^y V(t)\dd t-\int_{L-y}^L V(t)\dd t,
\]
and put $\omega(y)=V(y)-V(L-y)$, so that $\omega(L-y)=-\omega(y)$ and $A(y)=\int_0^y\omega$.

\begin{lemma}[First term]\label{lem:first}
For each fixed $y\in(0,L]$, as $k\to\infty$,
\[
w_-(y)=1+\frac1{2k}\int_0^yV+o(k^{-1}),\qquad w_+(y)=1+\frac1{2k}\int_{L-y}^LV+o(k^{-1}),
\]
hence
\begin{equation}\label{eq:first}
g(\lambda,y)=1+\frac{A(y)}{2k}+o(k^{-1}),\qquad A(y)=\int_0^yV-\int_{L-y}^LV .
\end{equation}
  
\end{lemma}
\begin{proof}
Split $1-e^{-2k(y-t)}$ in \eqref{eq:volterra}:
\[
w_-(y)=1-e^{-2ky}+\frac1{2k}\Bigl[\int_0^yV+\int_0^yV\,(w_--1)-\int_0^ye^{-2k(y-t)}V\,w_-\Bigr].
\]
By Lemma~\ref{lem:volterra}, for $k\ge1/2$, $|w_-|\le e^M$ and $|w_-(t)-1|\le e^{-2kt}+\frac{M}{2k}e^{M}$. Hence
\[
\begin{aligned}
\Bigl|\int_0^yV(w_--1)\Bigr|&\le\int_0^y|V(t)|e^{-2kt}\dd t+\frac{M^2e^{M}}{2k},\\
\Bigl|\int_0^ye^{-2k(y-t)}Vw_-\Bigr|&\le e^{M}\int_0^ye^{-2k(y-t)}|V(t)|\dd t.
\end{aligned}
\]
Both integrands on the right are dominated by $|V|\in L^1$ and tend to $0$ for almost every $t\in(0,y)$, so both integrals tend to $0$ by dominated convergence. Since $e^{-2ky}=o(k^{-1})$ for $y>0$, the formula for $w_-$ follows. For $w_+$ apply it to $\wt V$, noting $\int_0^yV(L-s)\dd s=\int_{L-y}^LV$. Finally $g=w_-/w_+$ and $1/(1+\varepsilon)=1-\varepsilon+O(\varepsilon^2)$ with $\varepsilon=O(k^{-1})$.
\end{proof}

\begin{theorem}[Asymptotic duality forces symmetry]\label{thm:asym}
Suppose that for each $y\in(0,L)$ there is $\kappa(y)>0$ with
\[
\frac{U_+(\lambda,y)}{U_-(\lambda,L-y)}=\kappa(y)\bigl(1+o(\lambda^{-1/2})\bigr)\qquad(\lambda\to\infty).
\]
Then $V=\wt V$ almost everywhere, $\kappa(y)=e^{(B(L)-B(y)+B(L-y))/2}$, and the duality is exact for every $\lambda>0$ (Theorem~\ref{thm:real}(i)). Moreover, for any $V$, the integrated asymmetry is read from the exit pair at high energy:
\[
A(y)=2\lim_{k\to\infty}k\,\bigl(g(\lambda,y)-1\bigr).
\]
  
\end{theorem}
\begin{proof}
By \eqref{eq:U}, $U_+(\lambda,y)/U_-(\lambda,L-y)=\kappa_0(y)\,g(\lambda,y)$ with $\kappa_0(y)=e^{(B(L)-B(y)+B(L-y))/2}$. Since $g\to1$ (Lemma~\ref{lem:first}), $\kappa=\kappa_0$, and then $g=1+o(k^{-1})$. Comparing with \eqref{eq:first}, $A(y)=0$ for every $y\in(0,L)$. $A$ is absolutely continuous with $A'(y)=V(y)-V(L-y)$ almost everywhere, so $V=\wt V$ almost everywhere. The limit formula is \eqref{eq:first}.
\end{proof}

\paragraph{The second high-energy term.} The next expansion has a rigorously controlled remainder provided $V$ is continuous at the four evaluation points.
\begin{proposition}[Second term]\label{prop:second}
Let $y\in(0,L)$ and suppose that $V$ has a representative continuous at $0$, $y$, $L-y$ and $L$ (one-sided at the endpoints). Then, as $\lambda=k^2\to\infty$,
\[
\log g(\lambda,y)=\frac{A(y)}{2k}-\frac{\omega(y)+\omega(0)}{4k^2}+o(k^{-2}),
\]
equivalently,
\begin{equation}\label{eq:second}
g(\lambda,y)=1+\frac{A(y)}{2k}+\frac1{k^2}\left[\frac{A(y)^2}{8}-\frac{\omega(0)+\omega(y)}{4}\right]+o(k^{-2}).
\end{equation}
\end{proposition}

\begin{proof}
Write $w=w_-$, $I(t)=\int_0^tV$, and put $\varepsilon(t)=w(t)-1+e^{-2kt}-I(t)/(2k)$. By Lemma~\ref{lem:volterra}, $|\varepsilon(t)|\le C/k$ for $k\ge1/2$ and all $t$, with $C$ depending only on $M$; by Lemma~\ref{lem:first} at the point $t$, $k\,\varepsilon(t)\to0$ for every $t\in(0,L]$. Hence, by dominated convergence, $\int_0^yV\varepsilon=o(1/k)$, and
\[
\int_0^yV(w-1)=-\int_0^ye^{-2kt}V(t)\dd t+\frac1{2k}\int_0^yV(t)I(t)\dd t+o(1/k)=-\frac{V(0)}{2k}+\frac{I(y)^2}{4k}+o(1/k),
\]
using $\int_0^yVI=I(y)^2/2$ and $\int_0^y2k\,e^{-2kt}V(t)\dd t\to V(0)$ (continuity of $V$ at $0$). In the same way, since $|w-1|\le e^{-2kt}+C/k$,
\[
\int_0^ye^{-2k(y-t)}V\,w=\int_0^ye^{-2k(y-t)}V(t)\dd t+O(k^{-2})+O(e^{-2ky})=\frac{V(y)}{2k}+o(1/k),
\]
by continuity of $V$ at $y$ ($V$ is bounded near $y$, which gives the $O(k^{-2})$). Inserting both in the split form of \eqref{eq:volterra} used in the proof of Lemma~\ref{lem:first},
\[
w_-(y)=1+\frac{I(y)}{2k}+\frac{I(y)^2}{8k^2}-\frac{V(0)+V(y)}{4k^2}+o(k^{-2}),\qquad
\log w_-(y)=\frac{I(y)}{2k}-\frac{V(0)+V(y)}{4k^2}+o(k^{-2}).
\]
For $w_+(y)$ apply this to $\wt V$, which is continuous at $0$ and $y$ because $V$ is continuous at $L$ and $L-y$: $\log w_+(y)=\frac1{2k}\int_{L-y}^LV-\frac{V(L)+V(L-y)}{4k^2}+o(k^{-2})$. Subtracting, and using $\log g=\log w_--\log w_+$, gives the result, since $V(y)-V(L-y)=\omega(y)$ and $V(0)-V(L)=\omega(0)$.
\end{proof}

\paragraph{All orders (formal, for smooth $V$).} Let $m_+=k+\sum_{n\ge1}a_nk^{-n}$ be the formal growing solution of $m'+m^2=k^2+V$, with $a_1=V/2$ and $a_n=-\frac12\bigl(a_{n-1}'+\sum_{i+j=n-1,\,i,j\ge1}a_ia_j\bigr)$; let $m_-(k)=m_+(-k)$ and $\Omega=\frac12(m_+-m_-)=k+\sum_{n\text{ odd}}a_nk^{-n}$. Then, up to $O(e^{-2ky})$,
\begin{equation}\label{eq:allorders}
\log g(\lambda,y)\ \sim\ \int_0^y\Omega-\int_{L-y}^L\Omega-\frac12\log\frac{\Omega(y)\,\Omega(0)}{\Omega(L-y)\,\Omega(L)} .
\end{equation}
The formal argument uses three standard identities: (a) $\log S_-(y)=\int_0^ym_+-\log\bigl(m_+(0)-m_-(0)\bigr)+O(e^{-2ky})$; (b) the even part satisfies $\frac12(m_++m_-)=-\frac12(\log\Omega)'$, by subtracting the Riccati equations for $m_\pm$; (c) $S_+(L-y)=S_-(y;\wt V)$ and $\Omega[\wt V](s)=\Omega[V](L-s)$, because the odd $a_n$ contain an even number of derivatives. With $\log(\Omega/k)=\sum_{n\text{ even}}l_nk^{-n}$, the coefficients of $\log g=\sum c_nk^{-n}$ are
\[
c_n=\int_0^ya_n-\int_{L-y}^La_n\ (n\text{ odd}),\qquad c_n=-\tfrac12\bigl[l_n(y)+l_n(0)-l_n(L-y)-l_n(L)\bigr]\ (n\text{ even}),
\]
with $a_3=(V''-V^2)/8$, $a_5=(2V^3-6VV''+V^{(4)}-5V'^2)/32$ (KdV densities), $l_2=V/2$, $l_4=V''/8-V^2/4$. In particular $c_1=\frac12\int_0^y\omega$ and $c_2=-\frac14(\omega(y)+\omega(0))$, in agreement with the proved expansion in Proposition~\ref{prop:second}; the linear part of $c_{2j}$ is $-2^{-2j}\bigl(\omega^{(2j-2)}(y)+\omega^{(2j-2)}(0)\bigr)$. Thus the odd orders integrate asymmetries of the KdV densities, whereas the even orders involve the even jets of the asymmetry at $0$ and at $y$.

\section{Neutral pairs: rigidity and nonsymmetric constructions}

\begin{figure}[t]
\centering
\begin{tikzpicture}[x=10.2cm,y=1cm,>=Stealth]
\draw[thick] (0,0)--(1,0);
\foreach \x in {0,.22,.5,.78,1}{\fill (\x,0) circle (1.5pt);}
\node[below=5pt] at (0,0) {$0$};
\node[below=5pt] at (.22,0) {$y$};
\node[below=5pt] at (.5,0) {$L/2$};
\node[below=5pt] at (.78,0) {$L-y$};
\node[below=5pt] at (1,0) {$L$};
\draw[->,thick] (.22,.2) .. controls (.36,.78) and (.64,.78) .. (.78,.2);
\node[above] at (.5,.71) {reflected reset pair};
\draw[dashed] (.5,-.25)--(.5,.64);
\end{tikzpicture}
\caption{The restart law supported at the reflected positions $y$ and $L-y$. Neutrality concerns the averaged exit measures, not independent copies with fixed starts.}
\label{fig:pairs}
\end{figure}
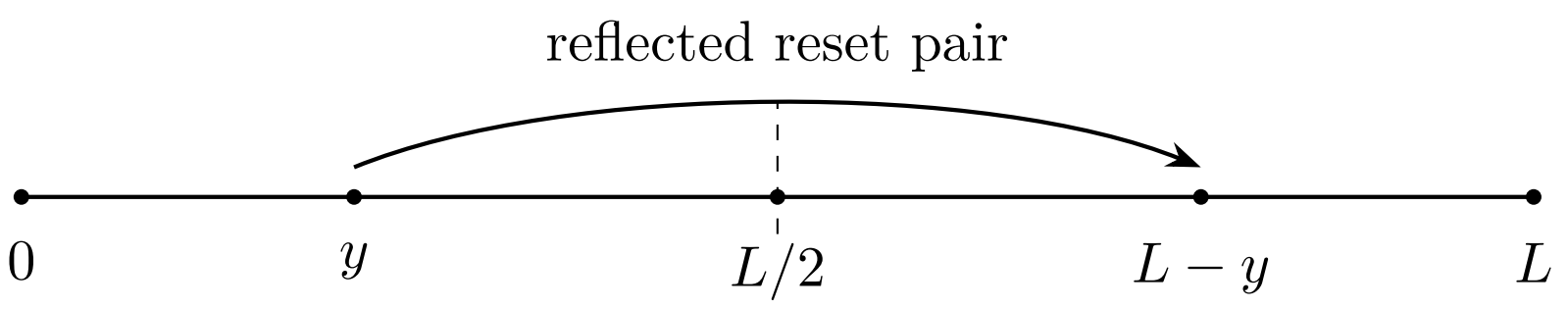

\begin{theorem}[Continuous realization of a spectral duality structure]\label{thm:real}
Define the duality factor and its orbit product by\[\kappa(y)=e^{(B(L)-B(y)+B(L-y))/2},\qquad K=e^{B(L)}.\]
\begin{enumerate}[label=(\roman*)]
\item If $V=\wt V$ almost everywhere, then for every $\lambda>0$ and $y\in(0,L)$
\[
U_+(\lambda,y)=\kappa(y)\,U_-(\lambda,L-y),\qquad \kappa(y)\kappa(L-y)=K,\qquad\kappa(L/2)=\sqrt K .
\]
\item Conversely, if there is a function $\kappa$ independent of $\lambda$ with $U_+(\lambda,y)=\kappa(y)U_-(\lambda,L-y)$ for all $y\in(0,L)$ and all $\lambda$ in an unbounded set, then $V=\wt V$ almost everywhere and $\kappa$ is as above.
\item If $V=\wt V$, then for every $y\in(0,L/2)$ the pair law with $w_2/w_1=\kappa(y)/\sqrt K$, equivalently $w_1:w_2=e^{B(y)/2}:e^{B(L-y)/2}$, is neutral, with $q_\lambda\equiv C^\ast=1/(1+\sqrt K)$.
\end{enumerate}
Only the Liouville potential has to be symmetric; the drift $\beta$ need not be. Given a symmetric $V$ and any positive solution $\psi_0$ of $\psi''=V\psi$, $\beta=2\psi_0'/\psi_0$ realizes it, and $e^{B/2}=\psi_0/\psi_0(0)$.
   
\end{theorem}
\begin{proof}
(i) By Lemma~\ref{lem:rev} with $\wt V=V$, $S_+(\lambda,L-y)=S_-(\lambda,y)$. By \eqref{eq:U},
\[
\frac{U_+(\lambda,y)}{U_-(\lambda,L-y)}=\frac{e^{[B(L)-B(y)]/2}S_-(\lambda,y)}{e^{-B(L-y)/2}S_+(\lambda,L-y)}=e^{(B(L)-B(y)+B(L-y))/2}=\kappa(y).
\]
Then $\kappa(y)\kappa(L-y)=e^{B(L)}$ and $\kappa(L/2)=e^{B(L)/2}$.

(ii) The hypothesis gives $g(\lambda,y)=\kappa(y)/\kappa_0(y)$ independent of $\lambda$ in the unbounded set, with $\kappa_0$ as in (i). Letting $\lambda\to\infty$ in that set, Lemma~\ref{lem:first} gives $g\equiv1$, so $\kappa=\kappa_0$. Now fix one $\lambda$ in the set: $\varphi=S_-(\lambda,\cdot)=S_+(\lambda,L-\cdot)$ solves both $\varphi''=(\lambda+V)\varphi$ and $\varphi''=(\lambda+\wt V)\varphi$, hence $(V-\wt V)\varphi=0$ almost everywhere, and $\varphi>0$ on $(0,L)$ (Lemma~\ref{lem:pos}). (Alternatively, the proof of Theorem~\ref{thm:asym} applies along the unbounded set.)

(iii) With $a_1=w_1e^{-B(y)/2}$ and $a_2=w_2e^{-B(L-y)/2}$ the stated weights give $a_1=a_2$. For $C=C^\ast$, $c=1$ in Proposition~\ref{prop:pairform}, and the condition reads $\Phi(y)+\Phi(L-y)=0$ with $\Phi=S_+-S_-$, which holds because $S_+(y)=S_-(L-y)$ and $S_+(L-y)=S_-(y)$ by Lemma~\ref{lem:rev}.

For the last assertion: $\beta=2\psi_0'/\psi_0$ gives $\beta'/2+\beta^2/4=\psi_0''/\psi_0=V$, and $B=2\log(\psi_0/\psi_0(0))$.
\end{proof}

\begin{theorem}[The weights of a neutral pair are forced]\label{thm:weights}
Let $0<y<L/2$. The pair $\{y,L-y\}$ carries a neutral pair law $w_1\delta_y+w_2\delta_{L-y}$, $w_1,w_2>0$, $w_1+w_2=1$, with some value $C\in(0,1)$, if and only if
\begin{equation}\label{eq:pairid}
S_-(\lambda,y)+S_-(\lambda,L-y)=S_+(\lambda,y)+S_+(\lambda,L-y)\qquad\text{for all }\lambda>0 .
\end{equation}
In that case the weights are unique, $w_1:w_2=e^{B(y)/2}:e^{B(L-y)/2}$, and $C=C^\ast$.
\end{theorem}
\begin{proof}
By Proposition~\ref{prop:pairform}, $a_1\bigl(S_+(y)-cS_-(y)\bigr)+a_2\bigl(S_+(L-y)-cS_-(L-y)\bigr)=0$ for all $\lambda>0$, with $c>0$. Use the Gelfand--Levitan (transformation operator) representation \cite{GelfandLevitan,Marchenko}: for $V\in L^1(0,L)$ there is a continuous kernel $K_V$ with
\[
S_-(\lambda,x;V)=\frac{\sinh kx}{k}+\int_0^xK_V(x,t)\frac{\sinh kt}{k}\dd t .
\]
By Lemma~\ref{lem:rev}, $S_+(\lambda,x;V)=S_-(\lambda,L-x;\wt V)$ has the same form, with the atom at $L-x$. So each of the four terms is $\int\frac{\sinh kt}{k}\,\rho(\dd t)$ for a finite measure $\rho$ on $[0,L]$: an atom at $y$ for $S_-(y)$ and $S_+(L-y)$, an atom at $L-y$ for $S_-(L-y)$ and $S_+(y)$, plus an absolutely continuous part. The identity therefore says $\int\sinh(kt)\,\rho(\dd t)=0$ for all $k>0$, where $\rho$ has atom $a_2-ca_1$ at $y$, atom $a_1-ca_2$ at $L-y$, and otherwise a density. Such a transform vanishes for all $k>0$ only if $\rho=0$ on $(0,L]$: extend $\rho$ to the odd measure $\tilde\rho$ on $[-L,L]$; then $\int e^{kt}\tilde\rho(\dd t)$ is entire in $k$ and vanishes for $k>0$, hence identically, and on the imaginary axis it is the Fourier transform of $\tilde\rho$. Since $y\ne L-y$, the two atoms vanish: $a_2=ca_1$ and $a_1=ca_2$, so $c^2=1$, $c=1$ and $a_1=a_2$. Then $C=C^\ast$, the weights are as stated, and the condition becomes \eqref{eq:pairid}. Conversely, with these weights $a_1=a_2$ and, for $C=C^\ast$, $c=1$, so \eqref{eq:pairid} is the condition of Proposition~\ref{prop:pairform}.
\end{proof}

\begin{theorem}[An interval of neutral pairs forces local symmetry]\label{thm:interval}
If \eqref{eq:pairid} holds for every $y$ in an open interval $J\subset(0,L/2)$, then $V(y)=V(L-y)$ for almost every $y\in J$. 
\end{theorem}
\begin{proof}
Let $G=S_--S_+$ and $h(y)=G(y)+G(L-y)$, for fixed $\lambda$. Both $G$ and $y\mapsto G(L-y)$ solve their equations in the Carathéodory sense, so $h'$ is absolutely continuous and, almost everywhere,
\[
h''(y)=(\lambda+V(y))G(y)+(\lambda+V(L-y))G(L-y).
\]
On $J$, $h\equiv0$, so $h''=0$ almost everywhere there and $G(L-y)=-G(y)$; hence $(V(y)-V(L-y))G(\lambda,y)=0$ for almost every $y\in J$. Take $\lambda_n=k_n^2\to\infty$; outside a null set the identity holds for all $n$. For fixed $y<L/2$, Lemma~\ref{lem:volterra} gives $2k\,e^{-k(L-y)}G(\lambda,y)\to-1$, so $G(\lambda_n,y)\ne0$ for large $n$, and $V(y)=V(L-y)$.
\end{proof}

\begin{theorem}[Pairs accumulating at an endpoint force global symmetry]\label{thm:endpoint}
If \eqref{eq:pairid} holds for $y=y_n$ with $y_n\to0$, then $V=\wt V$ almost everywhere on $[0,L]$. 
\end{theorem}
\begin{proof}
Fix $\lambda$ and let $h$ be as above; $h$ is $C^1$. Always $h(0)=G(0)+G(L)=-S_+(0)+S_-(L)=0$. Since $h(y_n)=0$, Rolle gives $\xi_n\in(0,y_n)$ with $h'(\xi_n)=0$, so $h'(0)=0$, that is $G'(0)=G'(L)$. Now $G'(0)=1-S_+'(0)$ and $G'(L)=S_-'(L)+1$, so $S_-'(L;V)=-S_+'(0;V)=S_-'(L;\wt V)$ by Lemma~\ref{lem:rev}. Also $S_-(L;V)=S_-(L;\wt V)$ by the same lemma. The two characteristic-function equalities imply the entire cross-product identity \eqref{eq:weyleq}. Lemma~\ref{lem:weylunique}, which reduces the problem to Borg's two-spectra theorem for real integrable potentials~\cite{Borg,LevitanGasymov} and~\cite[Theorem~4.4]{HrynivMykytyuk}, gives $V=\wt V$ almost everywhere.
\end{proof}

\begin{remark}[No global version near the midpoint]\label{rem:noglobal}
Theorem~\ref{thm:interval} gives only local symmetry. If $V=(H\,|\,H)$ (translated halves) and the half $H$ satisfies $H(x)=H(L/2-x)$ for $0\le x\le u_0$, with $u_0\le L/4$, then by Theorem~\ref{thm:families}(b) every pair with $y\in[L/2-u_0,L/2)$ is neutral, while $V$ need not be symmetric. 
\end{remark}

\begin{figure}[t]
\centering
\begin{tikzpicture}[x=.83cm,y=.65cm]
\foreach \i in {0,...,9}{\draw (\i,0) rectangle (\i+1,1);}
\node at (.5,.5) {$Z$};\node at (1.5,.5) {$Z'$};\node at (2.5,.5) {$Z$};\node at (3.5,.5) {$Z'$};\node at (4.5,.5) {$Z$};
\node at (5.5,.5) {$Z$};\node at (6.5,.5) {$Z'$};\node at (7.5,.5) {$Z$};\node at (8.5,.5) {$Z'$};\node at (9.5,.5) {$Z$};
\draw[line width=1.3pt] (5,-.18)--(5,1.2);
\draw[decorate,decoration={brace,amplitude=4pt}] (0,1.35)--(5,1.35) node[midway,above=5pt] {$H$};
\draw[decorate,decoration={brace,amplitude=4pt}] (5,1.35)--(10,1.35) node[midway,above=5pt] {$H$};
\end{tikzpicture}
\caption{A word representation of a periodic-half construction: if $L/2=2y+\rho$, then $H=(ZZ')^2Z$ with $|Z|=\rho$ and $|ZZ'|=y$. The full potential repeats the half, $V=(H\mid H)$. The words $Z,Z'$ need not be odd or antisymmetric.}
\label{fig:halves}
\end{figure}
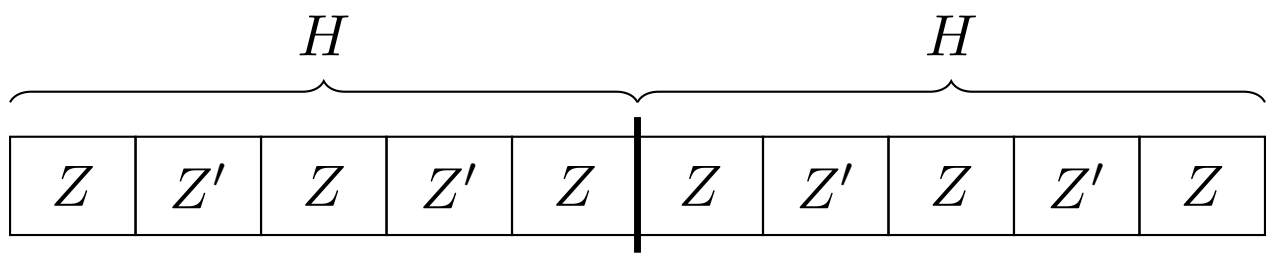

\begin{theorem}[Non-symmetric neutral pairs]\label{thm:families}
Each of the following gives a neutral pair, i.e.\ \eqref{eq:pairid}, for potentials that in general are not symmetric.
\begin{enumerate}[label=(\alph*)]
\item \emph{Commensurate periodic.} $V$ has period $p$, $L=Np$ and $y=mp$ with $1\le m<N/2$. (No condition $L/y\in\mathbb N$: e.g.\ $y=2L/9$.)
\item \emph{Translated halves.} $V=(H\,|\,H)$ with half $H=QPR$, $|Q|=|R|=u$, $y=L/2-u$ (so $L/4\le y<L/2$), and $Q$, $R$ with the same Dirichlet characteristic function, $[M_Q]_{12}\equiv[M_R]_{12}$ (e.g.\ $R=Q$ or $R=\wt Q$). Within the class of translated halves this condition is also necessary.
\item \emph{Translated halves with periodic half.} Any $y\in(0,L/2)$. $V=(H\,|\,H)$, and the half $H$ has period $y$: $H(x+y)=H(x)$ for $0\le x\le L/2-y$.
\end{enumerate}
Since \eqref{eq:pairid} is an identity between entire functions of $\lambda$, and $S_\pm(\lambda;V+\gamma)=S_\pm(\lambda+\gamma;V)$, each example stays neutral after adding a constant to $V$. Adding a sufficiently large constant makes a bounded $V$ positive, and a positive potential has a positive solution of $\psi''=V\psi$ (the solution with $\psi(0)=1$, $\psi'(0)=0$ is convex and increasing), so Theorem~\ref{thm:real} realizes it by a drift. In particular, taking bounded examples in each family and choosing a sufficiently large constant gives positive realizations by diffusions.
  
\end{theorem}
\begin{proof}
In each case we compute the four terms of \eqref{eq:pairid}. Recall $S_-(x;V)$ is the $(1,2)$ entry of the transfer matrix of $V|_{[0,x]}$, and $S_+(x;V)=S_-(L-x;\wt V)$ (Lemma~\ref{lem:rev}). By Lemma~\ref{lem:revmat}, $S_-(x;\wt V)$ is the $(1,2)$ entry of the transfer matrix of the \emph{reversal} of $\wt V|_{[0,x]}$, which is $V|_{[L-x,L]}$. So
\begin{equation}\label{eq:fourterms}
S_+(y)=[\,M_{V|[y,L]}\,]_{12},\qquad S_+(L-y)=[\,M_{V|[L-y,L]}\,]_{12}.
\end{equation}
(These are just the Dirichlet characteristic functions of the four subintervals $[0,y]$, $[0,L-y]$, $[y,L]$, $[L-y,L]$.)

(a) Translation by a multiple of $p$ does not change the potential, so $V|_{[y,L]}$ is a copy of $V|_{[0,L-y]}$ and $V|_{[L-y,L]}$ a copy of $V|_{[0,y]}$. By \eqref{eq:fourterms}, $S_+(y)=S_-(L-y)$ and $S_+(L-y)=S_-(y)$.

(b) Write $P,Q,R$ for the transfer matrices of the pieces. Then $[0,y]=QP$, $[0,L-y]=QPRQ$, $[y,L]=RQPR$, $[L-y,L]=PR$, so with the product rule
\[
\text{\eqref{eq:pairid} residual}=[PQ]_{12}+[QRPQ]_{12}-[RPQR]_{12}-[RP]_{12}.
\]
Let $X=RPQ$ and $t=\tr X=\tr(PQR)$. By Cayley--Hamilton, $X=tI-X^{-1}$, so $QX=tQ-QX^{-1}=tQ-(RP)^{-1}$ and $XR=tR-X^{-1}R=tR-(PQ)^{-1}$. Hence
\[
\begin{aligned}
&PQ+QRPQ-RP-RPQR\\
&\quad=\bigl(PQ+(PQ)^{-1}\bigr)-\bigl(RP+(RP)^{-1}\bigr)+t(Q-R)\\
&\quad=\bigl(\tr(PQ)-\tr(RP)\bigr)I+t(Q-R),
\end{aligned}
\]
and the $(1,2)$ entry is
\begin{equation}\label{eq:halves}
\text{residual}=\bigl([Q]_{12}-[R]_{12}\bigr)\,\tr(M_H).
\end{equation}
So $[Q]_{12}\equiv[R]_{12}$ is sufficient. Since $\tr M_H\sim e^{kL/2}$ (Lemma~\ref{lem:volterra}), it is not identically zero, and the entire function $[Q]_{12}-[R]_{12}$ must vanish: the condition is necessary within the class. $R=\wt Q$ works by Lemma~\ref{lem:revmat}.

(c) Write $L/2=qy+\rho$ with $q\ge1$, $0\le\rho<y$, and let $Z=H|_{[0,\rho]}$, $Z'=H|_{[\rho,y]}$, $a=M_Z$, $b=M_{Z'}$, $u=ba$. Periodicity gives $H=(ZZ')^qZ=Z(Z'Z)^q$, so, as a word, $V=(ZZ')^qZ\,Z(Z'Z)^q$. Hence $V|_{[0,y]}=ZZ'$ and $V|_{[0,L-y]}=(ZZ')^qZZ(Z'Z)^{q-1}$, with matrices
\[
ba\qquad\text{and}\qquad (ab)^{q-1}a\,a\,u^q=au^{q-1}au^q,
\]
using $(ab)^ja=a(ba)^j=au^j$. Reversing the word, $\wt V=(\wt Z\wt Z')^q\wt Z\,\wt Z(\wt Z'\wt Z)^q$ is the same word with $Z,Z'$ replaced by $\wt Z,\wt Z'$, whose matrices are $Ea^{-1}E$, $Eb^{-1}E$ (Lemma~\ref{lem:revmat}). Substituting in the two matrices above gives $E(ab)^{-1}E$ and $E\bigl[(ab)^qa(ab)^{q-1}a\bigr]^{-1}E$, and for $\det X=1$ the $(1,2)$ entry of $EX^{-1}E$ is $X_{12}$. With $S_+(x)=S_-(L-x;\wt V)$, the four terms are
\[
S_-(y)=[ba]_{12},\quad S_-(L-y)=[au^{q-1}au^q]_{12},\quad S_+(L-y)=[ab]_{12},\quad S_+(y)=[au^qau^{q-1}]_{12},
\]
and \eqref{eq:pairid} is $[au^qau^{q-1}-au^{q-1}au^q]_{12}=[ba-ab]_{12}$, which follows from Lemma~\ref{lem:matrix}.
\end{proof}

\begin{lemma}[A matrix identity]\label{lem:matrix}
For $2\times2$ matrices $a,b$ with determinant $1$, $u=ba$ and $q\ge1$:
\[
au^qau^{q-1}-au^{q-1}au^q=ba-ab.
\]
\end{lemma}
\begin{proof}
Put $F_j=au^{j+1}au^j-au^jau^{j+1}$ for $j\ge0$. Since $\det u=1$, Cayley--Hamilton gives
$u^{j+2}=(\tr u)u^{j+1}-u^j$. Substituting this recurrence in $F_{j+1}$, the terms involving $\tr u$ cancel and $F_{j+1}=F_j$. Thus $F_{q-1}=F_0$. Now
$F_0=abaa-aaba=ba-ab$, where the last equality follows from $a^2=(\tr a)a-I$.
\end{proof}

\begin{remark}[The simplest non-symmetric neutral midpoint]\label{rem:midpoint}
The limiting case $u=0$ of Theorem~\ref{thm:families}(b), $V(y+L/2)=V(y)$, concerns the midpoint rather than a pair: the midpoint is neutral by Theorem~\ref{thm:midpoint}(b) and Lemma~\ref{lem:revmat}, e.g.\ for $V=(v_1,v_2\,|\,v_1,v_2)$ piecewise constant with $v_1\ne v_2$. The explicit Weyl contrast appears in Example~\ref{ex:weyl}. 
\end{remark}

\begin{corollary}[A second necessary condition]\label{cor:necessary}
If the pair $\{y,L-y\}$, $0<y<L/2$, is neutral, then $\int_0^y\omega=0$. If moreover $V$ has a representative continuous at $0$, $y$, $L-y$ and $L$, then $\omega(y)=\omega(0)$, that is, $V(y)-V(L-y)=V(0)-V(L)$.
\end{corollary}

\begin{proof}
By Theorem~\ref{thm:weights}, neutrality gives \eqref{eq:pairid}, so $S_-(L-y)-S_+(y)=S_+(L-y)-S_-(y)$. By Lemma~\ref{lem:volterra} the right-hand side is $O(e^{ky}/k)$, while $S_+(y)\ge e^{k(L-y)}/(4k)$ for large $k$; hence $g(\lambda,L-y)=S_-(L-y)/S_+(y)=1+O(e^{-k(L-2y)})=1+o(k^{-2})$. By Lemma~\ref{lem:first} at the point $L-y$, $A(L-y)=0$, and $A(L-y)=\int_0^{L-y}\omega=\int_0^y\omega$ because $\omega(L-x)=-\omega(x)$. By Proposition~\ref{prop:second} at the point $L-y$ (its hypotheses are symmetric in $y\leftrightarrow L-y$), $\omega(L-y)+\omega(0)=0$, and $\omega(L-y)=-\omega(y)$.
\end{proof}


\section{Midpoint Weyl rigidity and the discrete analogue}\label{sec:weyl}
\subsection{Midpoint neutrality and the Weyl response}\label{subsec:midweyl}
For a single reset point, the high-energy behaviour of the Dirichlet functions restricts all-rate neutrality to the midpoint $c=L/2$. Equal half-interval Dirichlet spectra need not imply equal Weyl functions.
\begin{theorem}[Individual midpoint neutrality]\label{thm:midpoint}
Let $y\in(0,L)$. The following hold.
\begin{enumerate}[label=(\alph*)]
\item If $q_\lambda(y)$ is independent of $\lambda>0$, then $y=c:=L/2$.
\item The point $c$ is neutral, that is $q_\lambda(c)$ is independent of $\lambda>0$, if and only if
\begin{equation}\label{eq:midchar}
S_-(\lambda,c;V)=S_-(\lambda,c;\wt V)\qquad\text{for all }\lambda\in\mathbb C,
\end{equation}
that is, the Dirichlet characteristic functions of $V|_{[0,c]}$ and of $\wt V|_{[0,c]}$ (the right half reflected) coincide; and this holds if and only if the two halves have the same Dirichlet spectrum on $[0,c]$.
\item In that case $q_\lambda(c)=C^\ast=1/(1+e^{B(L)/2})$ for every $\lambda>0$.
\end{enumerate}
By Proposition~\ref{prop:single}, neutrality for every independent renewal clock is equivalent to the condition in (a)--(b).
\end{theorem}

\begin{proof}
By \eqref{eq:U}, $U_+(\lambda,y)/U_-(\lambda,y)=e^{B(L)/2}\rho(\lambda)$ with $\rho(\lambda)=S_-(\lambda,y)/S_+(\lambda,y)$, and $q_\lambda(y)=1/(1+e^{B(L)/2}\rho(\lambda))$. Since $U_\pm>0$, $q_\lambda(y)$ is constant if and only if $\rho$ is constant.

(a) By Lemma~\ref{lem:rev}, $S_+(\lambda,y;V)=S_-(\lambda,L-y;\wt V)$. By Lemma~\ref{lem:order} applied to $V$ on $[0,y]$ and to $\wt V$ on $[0,L-y]$,
\[
\rho(\lambda)=e^{k(2y-L)}\,\bigl(1+o(1)\bigr)\qquad(\lambda\to+\infty).
\]
If $y<L/2$ this tends to $0$, and if $y>L/2$ to $+\infty$; neither is compatible with a constant $\rho>0$. Hence $y=c$.

(b) At $y=c$, $\rho(\lambda)=S_-(\lambda,c;V)/S_-(\lambda,c;\wt V)\to1$ by Lemma~\ref{lem:order}, so $\rho$ is constant if and only if $\rho\equiv1$ on $(0,\infty)$, which is \eqref{eq:midchar} for $\lambda>0$; both sides are entire, so \eqref{eq:midchar} then holds on $\mathbb C$. $S_-(\lambda,c;V)$ depends only on $V|_{[0,c]}$ and $S_-(\lambda,c;\wt V)$ only on $\wt V|_{[0,c]}$. Equal functions have equal zeros. Conversely, if the two Dirichlet spectra coincide (with multiplicities; Dirichlet eigenvalues are simple), the two functions have the same zeros, and the Hadamard argument in the proof of Lemma~\ref{lem:weylunique} gives \eqref{eq:midchar}.

(c) With $\rho\equiv1$, $q_\lambda(c)=1/(1+e^{B(L)/2})$.
\end{proof}

For the outward logarithmic derivatives $m_L=S_-'/S_-$ and $m_R=-S_+'/S_+$, with derivatives in the intrinsic coordinate, put
\begin{equation}\label{eq:weyl}
\partial_y\log U_+=m_L-\frac{\beta}{2},\qquad
\partial_y\log U_-=-m_R-\frac{\beta}{2},\qquad
\Gamma:=\partial_y\log(U_-U_+)+\beta=m_L-m_R.
\end{equation}
By Lemma~\ref{lem:rev}, the two logarithmic derivatives are the Weyl functions of the outward halves, expressed in the same inward coordinate:
\begin{equation}\label{eq:mR}
m_L(\lambda;y)=\frac{s'(\lambda,y;V)}{s(\lambda,y;V)},\qquad
m_R(\lambda;y)=\frac{s'(\lambda,L-y;\wt V)}{s(\lambda,L-y;\wt V)}.
\end{equation}
\begin{theorem}[Weyl reflection rigidity]\label{thm:weyl}
If $\Gamma(\lambda;c)=0$ for all $\lambda$ in a subset of $(0,\infty)$ with an accumulation point in $(0,\infty)$ (for instance a nonempty open interval), then $V=\wt V$ almost everywhere on $(0,L)$. Conversely, if $V=\wt V$ almost everywhere, then $\Gamma(\lambda;c)=0$ for every $\lambda>0$. In particular $\Gamma(\cdot;c)\equiv0$ is equivalent to midpoint neutrality \eqref{eq:midchar} together with equality of the Dirichlet--Neumann spectra of the two halves. The statement concerns the Liouville potential $V$; the coefficients $(D,b)$ need not be symmetric.
\end{theorem}

\begin{proof}
By \eqref{eq:mR} at $y=c$ (so $L-y=c$), $\Gamma(\lambda;c)=0$ is \eqref{eq:weyleq} with $W_1=V|_{[0,c]}$ and $W_2=\wt V|_{[0,c]}$; for $\lambda>0$ the denominators are positive (Lemma~\ref{lem:pos}). Lemma~\ref{lem:weylunique} gives $V(y)=\wt V(y)=V(L-y)$ for almost every $y\in(0,c)$, hence almost everywhere on $(0,L)$. Conversely, if $V=\wt V$, then $W_1=W_2$ and the two Weyl functions coincide. The last assertion is the first part of Lemma~\ref{lem:weylunique}: $\Gamma\equiv0$ if and only if $s_1=s_2$ and $s_1'=s_2'$; and $s_1'=s_2'$ is equivalent to equal Dirichlet--Neumann spectra by the Hadamard argument of that lemma applied to $s'$, with the normalization $2e^{-kc}s'\to1$ of Lemma~\ref{lem:order}.
\end{proof}

\begin{example}[Midpoint neutral but not Weyl rigid]\label{ex:weyl}
Let $V=(H\,|\,H)$ with $H=(v_1,v_2)$ piecewise constant, of lengths $l_1,l_2$, $v_1,v_2>0$ with $v_1\ne v_2$, and $k_1=\sqrt{\lambda+v_1}$, $k_2=\sqrt{\lambda+v_2}$. The midpoint is neutral (Remark~\ref{rem:midpoint}, Theorem~\ref{thm:midpoint}), but
\[
\Gamma(\lambda;c)=\frac{(v_2-v_1)\sinh(k_1l_1)\sinh(k_2l_2)}{k_1k_2\,S_-(\lambda,c)}\neq0\qquad(\lambda>0),
\]
and $S_-'(\lambda,c;V)-S_-'(\lambda,c;\wt V)=(v_2-v_1)\sinh(k_1l_1)\sinh(k_2l_2)/(k_1k_2)$: the two halves have the same Dirichlet spectrum and different Dirichlet--Neumann spectra. Indeed, with $M_H=M_{v_2}M_{v_1}$, the half $\wt V|_{[0,c]}$ is $H$ reversed, so by Lemma~\ref{lem:revmat} its transfer matrix has the diagonal entries of $M_H$ exchanged; hence $\Gamma=([M_H]_{22}-[M_H]_{11})/[M_H]_{12}$, and a direct product of the two constant-piece matrices gives $[M_H]_{22}-[M_H]_{11}=\sinh(k_1l_1)\sinh(k_2l_2)\,(k_2/k_1-k_1/k_2)$. For positive $v_1,v_2$ the potential is realized by a diffusion through a positive solution of $\psi''=V\psi$; an arbitrary bounded-below example can first be shifted by a sufficiently large constant without changing its neutrality.
\end{example}

\begin{proposition}[The first spatial jet identifies the intrinsic coefficients]\label{prop:jet}
Fix $y_0\in(0,L)$, $L$ and the intrinsic coordinate. Suppose that the functions $\lambda\mapsto U_\pm(\lambda,y_0)$ and $\lambda\mapsto\partial_yU_\pm(\lambda,y_0)$ are known for $\lambda$ in a nonempty open interval $I\subset(0,\infty)$. Then $V$ is determined almost everywhere on $(0,L)$ and $\beta$ is determined on $[0,L]$.
\end{proposition}

\begin{proof}
\emph{Continuation.} By \eqref{eq:U}, $U_+(\lambda,y_0)=e^{[B(L)-B(y_0)]/2}S_-(\lambda,y_0)/\Delta(\lambda)$, and similarly for $U_-$ and the $y$-derivatives. These are quotients of entire functions of $\lambda$ with $\Delta>0$ on $(0,\infty)$, hence real-analytic on $(0,\infty)$; their values on $I$ determine them on $(0,\infty)$ by the identity theorem. So the logarithmic derivatives
\[
\ell_\pm(\lambda)=\frac{\partial_yU_\pm(\lambda,y_0)}{U_\pm(\lambda,y_0)}
\]
are known for every $\lambda>0$. By \eqref{eq:weyl}, $\ell_+=m_L-\beta(y_0)/2$ and $\ell_-=-m_R-\beta(y_0)/2$.

\emph{The value $\beta(y_0)$.} By Lemma~\ref{lem:order} and \eqref{eq:mR}, $m_L(\lambda;y_0)=k+o(1)$ as $\lambda=k^2\to\infty$. Hence
\[
\beta(y_0)=2\lim_{k\to\infty}\bigl(k-\ell_+(k^2)\bigr)
\]
is determined, and then so are $m_L=\ell_++\beta(y_0)/2$ and $m_R=-\ell_--\beta(y_0)/2$ on $(0,\infty)$.

\emph{The potential.} If two admissible drifts $\beta,\hat\beta$ give the same data, their potentials $V,\hat V$ have the same $m_L(\cdot;y_0)$, which is \eqref{eq:weyleq} for $V|_{[0,y_0]}$ and $\hat V|_{[0,y_0]}$ (with $c=y_0$); Lemma~\ref{lem:weylunique} gives $V=\hat V$ on $(0,y_0)$. Likewise equality of $m_R(\cdot;y_0)$ gives $\wt V=\wt{\hat V}$ on $(0,L-y_0)$, that is $V=\hat V$ on $(y_0,L)$.

\emph{The drift.} $\beta$ and $\hat\beta$ belong to $W^{1,1}(0,L)$, satisfy $\beta'=2V-\beta^2/2$ almost everywhere (the same $V$) and $\beta(y_0)=\hat\beta(y_0)$. The difference $\delta=\beta-\hat\beta$ is absolutely continuous with $\delta'=-\tfrac12(\beta+\hat\beta)\delta$, and $\beta+\hat\beta$ is bounded, so Gronwall's inequality from $y_0$ in both directions gives $\delta\equiv0$.
\end{proof}

Exact-data recovery is stronger than midpoint isospectrality. Spectral deformations give isospectral families in other settings~\cite{GesztesySimonTeschl}; they are not needed for the separation in Example~\ref{ex:weyl}.

\subsection{A finite birth--death chain}
For a nearest-neighbour chain on $\{0,1,\dots,a\}$ with absorbing endpoints, let $p_j\in(0,1)$ be the probability of stepping right from site $j$ and $q_j=1-p_j$. Let $P$ be the killed transition matrix, $\Delta(s)=\det(I-sP)$, $d_L(j;s)=\det(I-sP_{[1,j]})$ and $d_R(j;s)=\det(I-sP_{[j,a-1]})$, with empty determinants one. The discounted exit channels are $U_+(s,z)=\mathbb E_z[s^T;X_T=a]$ and $U_-(s,z)=\mathbb E_z[s^T;X_T=0]$, for $0<s<1$. The Jacobi-matrix approach is classical~\cite{GantmacherKrein}.
\begin{lemma}[Renewal clocks for the chain]\label{lem:discclock}
For an interior site $z$, the eventual exit probability at $0$ is the same for every independent discrete renewal clock if and only if $U_+(s,z)/U_-(s,z)$ is constant on $(0,1)$.
\end{lemma}

\begin{proof}
The renewal argument of Proposition~\ref{prop:single} gives, for a clock whose inter-reset time $\tau$ has survival function $W(n)=\mathbb P(\tau>n)$ on $\{0,1,2,\dots\}$, a reset at time $\tau$ taking place before the $\tau$-th step (so that a reset and an absorption never coincide), $Q_z(W)=\sum_nW(n)\nu(n)/\sum_nW(n)(\nu+\eta)(n)$ with $\nu(n)=\mathbb P_z(T=n,X_T=0)$ and $\eta(n)=\mathbb P_z(T=n,X_T=a)$. The geometric clocks $W(n)=s^n$ give $U_\pm(s,z)$. If $U_+/U_-=\kappa$ on $(0,1)$, the power series $\sum_n(\eta(n)-\kappa\nu(n))s^n$ vanishes on $(0,1)$, so $\eta=\kappa\nu$, and then $Q_z(W)=1/(1+\kappa)$ for every $W$. The converse is the case of geometric clocks.
\end{proof}

\begin{lemma}[Jacobi determinant formulas]\label{lem:jacobi}
For $z\in\{1,\dots,a-1\}$ and $s\in(0,1)$,
\begin{equation}\label{eq:jacobi}
U_+(s,z)=\frac{s^{a-z}\bigl(\prod_{j=z}^{a-1}p_j\bigr)\,d_L(z-1;s)}{\Delta(s)},\qquad
U_-(s,z)=\frac{s^{z}\bigl(\prod_{j=1}^{z}q_j\bigr)\,d_R(z+1;s)}{\Delta(s)} .
\end{equation}
\end{lemma}

\begin{proof}
First-step analysis gives $(I-sP)u_+=s\,p_{a-1}e_{a-1}$ for the vector $u_+=(U_+(s,z))_z$; $I-sP$ is invertible because $P$ is substochastic with spectral radius less than one. By Cramer's rule, $U_+(s,z)=\det A_z/\Delta(s)$, where $A_z$ is $I-sP$ with column $z$ replaced by $s\,p_{a-1}e_{a-1}$. Expanding along that column,
$\det A_z=s\,p_{a-1}(-1)^{(a-1)+z}\det N$, where $N$ is $I-sP$ without row $a-1$ and column $z$. The rows $1,\dots,z-1$ of $N$ have entries only in the columns $1,\dots,z-1$, where they form $I-sP_{[1,z-1]}$; in the rows $z,\dots,a-2$ and the columns $z+1,\dots,a-1$ the matrix is lower triangular with diagonal entries $(I-sP)_{j,j+1}=-s\,p_j$. Hence $N$ is block lower triangular and
\[
\det N=d_L(z-1;s)\prod_{j=z}^{a-2}(-s\,p_j),\qquad
\det A_z=(-1)^{(a-1)+z}(-1)^{a-1-z}\,s^{a-z}\Bigl(\prod_{j=z}^{a-1}p_j\Bigr)d_L(z-1;s),
\]
and the sign is $(-1)^{2(a-1)}=1$. The formula for $U_-$ follows by reflecting the chain ($j\mapsto a-j$, $p\leftrightarrow q$).
\end{proof}

\begin{theorem}[Discrete midpoint isospectrality]\label{thm:disc}
An interior site $z$ is neutral for every independent discrete renewal clock if and only if $a=2z$ and $d_L(z-1;s)=d_R(z+1;s)$ as polynomials in $s$. In that case the eventual left-exit probability is
\[
C^\ast=\frac{1}{1+\kappa_c},\qquad \kappa_c=\frac{\prod_{j=z}^{a-1}p_j}{\prod_{j=1}^{z}q_j} ,
\]
which for constant $p_j=p$ is $1/(1+(p/q)^{a/2})$.
\end{theorem}

\begin{proof}
By Lemma~\ref{lem:jacobi}, $U_+(s,z)/U_-(s,z)=\kappa_c\,s^{a-2z}\,d_L(z-1;s)/d_R(z+1;s)$, and both determinants equal $1$ at $s=0$. By Lemma~\ref{lem:discclock}, neutrality means that this ratio is constant on $(0,1)$. If $a\ne2z$, it tends to $0$ or to $+\infty$ as $s\downarrow0$, so it is not constant. If $a=2z$, it tends to $\kappa_c$, so it is constant if and only if $d_L(z-1;s)=d_R(z+1;s)$ on $(0,1)$, that is, as polynomials. The value is $1/(1+\kappa_c)$.
\end{proof}

\begin{example}[$a=8$]
For a $3\times3$ tridiagonal matrix with unit diagonal, $\det=1-x_1y_1-x_2y_2$ (product of the off-diagonal pairs). So $d_L(3;s)=1-s^2(p_1q_2+p_2q_3)$ and $d_R(5;s)=1-s^2(p_5q_6+p_6q_7)$, and the midpoint $z=4$ is neutral if and only if $p_1q_2+p_2q_3=p_5q_6+p_6q_7$. An explicit non-symmetric solution (symmetry would mean $p_j=q_{8-j}$) is
\[
(p_1,\dots,p_7)=\bigl(\tfrac12,\tfrac12,\tfrac12,\tfrac12,\tfrac25,\tfrac27,\tfrac14\bigr),\qquad p_1q_2+p_2q_3=\tfrac12=\tfrac25\cdot\tfrac57+\tfrac27\cdot\tfrac34 .
\]
\end{example}

\begin{proposition}[Recovery from adjacent exit pairs]\label{prop:three}
Let $2\le z\le a-2$. The functions $U_+(\cdot,z)$, $U_+(\cdot,z+1)$, $U_-(\cdot,z-1)$ and $U_-(\cdot,z)$ on $(0,1)$ determine all the transition probabilities $p_1,\dots,p_{a-1}$.
\end{proposition}

\begin{proof}
Put $e_j=p_{j-1}q_j$ ($2\le j\le a-1$). Expanding along the last row, $d_L(j)=d_L(j-1)-s^2e_jd_L(j-2)$ for $j\ge2$, with $d_L(0)=d_L(1)=1$. Hence the rational functions $r_j=d_L(j)/d_L(j-1)$ satisfy
\begin{equation}\label{eq:cf}
r_1=1,\qquad r_j=1-\frac{s^2e_j}{r_{j-1}}\quad(j\ge2),\qquad r_j(0)=1 .
\end{equation}
By Lemma~\ref{lem:jacobi}, $R_L(s):=U_+(s,z+1)/U_+(s,z)=r_z(s)/(s\,p_z)$. Since $r_z(0)=1$, $p_z=\lim_{s\downarrow0}1/(sR_L(s))$, and then $r_z=s\,p_zR_L$ is known. By \eqref{eq:cf}, $e_j=\lim_{s\downarrow0}(1-r_j(s))/s^2$ and $r_{j-1}=s^2e_j/(1-r_j)$; descending from $j=z$ to $j=2$ we obtain $e_z,\dots,e_2$. Now $q_z=1-p_z$ is known, and $p_{j-1}=e_j/q_j$, $q_{j-1}=1-p_{j-1}$ give $p_{z-1},\dots,p_1$.

Symmetrically, $d_R(j)=d_R(j+1)-s^2f_jd_R(j+2)$ with $f_j=p_jq_{j+1}$, and $U_-(s,z-1)/U_-(s,z)=\bigl(d_R(z)/d_R(z+1)\bigr)/(s\,q_z)$, which by the same descent gives $f_z,\dots,f_{a-2}$; then $q_{j+1}=f_j/p_j$, $p_{j+1}=1-q_{j+1}$ give $p_{z+1},\dots,p_{a-1}$. (The limit $\lim sU_-(s,z-1)/U_-(s,z)$ gives $1/q_z$, a consistency check.)
\end{proof}

The data are four full functions of the discount parameter, not four scalar measurements.

\section{Two-rate reconstruction in physical coordinates}\label{sec:inverse}
In physical coordinates let $u_\pm(\lambda,x)$ solve $D u''+b u'=\lambda u$ with boundary values $(1,0)$ for $u_-$ and $(0,1)$ for $u_+$. Set $q_\lambda=u_-/(u_-+u_+)$; then $q_\lambda(0)=1$, $q_\lambda(\ell)=0$ and $q_\lambda'<0$. To see strict monotonicity, the Wronskian $\mathcal W=u_-u_+'-u_-'u_+$ satisfies $\mathcal W'=-g\mathcal W$, with $g=b/D$, and $\mathcal W(0)=u_+'(0)>0$. With $F=u_-+u_+$, $q_\lambda'=-\mathcal W/F^2<0$ on $[0,\ell]$. 
For a $C^3$ function $f$ with $f'\ne0$, its Schwarzian derivative is
\[
\mathcal S(f)=\frac{f'''}{f'}-\frac32\left(\frac{f''}{f'}\right)^2.
\]

\subsection{Classical Schwarzian identity and data consistency}

\begin{proposition}[Schwarzian of the reset exit profile]\label{prop:schwarz}
In the physical coordinate, with $g=b/D$ and $q_\lambda=u_-/(u_-+u_+)$, for every $\lambda\ge0$
\[
\cS(q_\lambda)=-\frac{2\lambda}D-g'-\frac12g^2 .
\]
Consequently, for $\lambda_1\ne\lambda_2$, $1/D=\bigl(\cS(q_{\lambda_1})-\cS(q_{\lambda_2})\bigr)/\bigl(2(\lambda_2-\lambda_1)\bigr)$, and exact data from a diffusion must satisfy
\[
\frac{\cS(q_{\lambda_1})-\cS(q_{\lambda_2})}{\lambda_2-\lambda_1}>0\quad\text{on }[0,\ell].
\]
For $\lambda=0$, $u_-+u_+=1$ and $g=-q_0''/q_0'$.   
\end{proposition}
\begin{proof}
$u_\pm$ solve $y''+gy'-(\lambda/D)y=0$, and so does $F=u_-+u_+$; $q=u_-/F$ is a ratio of two solutions. For $y''+Py'+Qy=0$ with solutions $y_1,y_2$ and $q=y_1/y_2$: the Wronskian $W=y_1'y_2-y_1y_2'$ satisfies $W'=-PW$, and $q'=W/y_2^2$. So $q''/q'=-P-2y_2'/y_2$, and differentiating once more and using the equation for $y_2$,
\[
\cS(q)=\Bigl(\frac{q''}{q'}\Bigr)'-\frac12\Bigl(\frac{q''}{q'}\Bigr)^2=2Q-P'-\frac12P^2 .
\]
With $P=g$, $Q=-\lambda/D$ this is the stated identity. Subtracting at two rates eliminates $g$. At $\lambda=0$, $F\equiv1$ because $u_-+u_+$ is harmonic with boundary values $1$, so $q_0=u_-$ and $q_0''+gq_0'=0$.
\end{proof}


\subsection{Two exact profiles identify both coefficients}
\begin{theorem}[Exact identification by two distinct rates]\label{thm:two}
Let $D\in C^2[0,\ell]$ be strictly positive and $b\in C^1[0,\ell]$. Let $\lambda_1\ne\lambda_2$ be nonnegative rates and suppose the exact full profiles $q_{\lambda_1},q_{\lambda_2}$ arise from the same diffusion. Writing $g=b/D$, the diffusivity is recovered pointwise from
\begin{equation}\label{eq:Drec}
 \frac1{D(x)}=\frac{\mathcal S(q_{\lambda_1})(x)-\mathcal S(q_{\lambda_2})(x)}{2(\lambda_2-\lambda_1)}.
\end{equation}
With $q=q_{\lambda_1}$ and this $D$, the boundary problem
\begin{equation}\label{eq:wbvp}
 w''-\frac{q''}{q'}w'+\frac{\lambda_1}{D}w=0,\qquad w(0)=w(\ell)=1
\end{equation}
has exactly one solution, which is positive, and
\begin{equation}\label{eq:brec}
 b=D\left(-\frac{q''}{q'}+2\frac{w'}w\right).
\end{equation}
Thus $(D,b)$ is uniquely identified in the fixed physical coordinate.
\end{theorem}
\begin{proof}
The Schwarzian identity in Proposition~\ref{prop:schwarz} gives \eqref{eq:Drec}. Set $F=u_-+u_+$ and $w=1/F$. Since $u_-=qF$, subtracting the equations for $u_-$ and $F$ yields $q''+(2F'/F+g)q'=0$, hence \eqref{eq:brec}. Substitution of $w=1/F$ into the equation for $F$ gives \eqref{eq:wbvp}; endpoint values are one because $F=1$ at both endpoints.

Existence and positivity follow from $F(x)=\mathbb E_x[e^{-\lambda_1T}]>0$. To prove uniqueness, set $p=-1/q'>0$ and $r=\lambda_1/D\ge0$. Then \eqref{eq:wbvp} is $(pw')'+rp w=0$. For every smooth $f$ vanishing at both endpoints, the ground-state identity is
\begin{equation}\label{eq:ground}
\int_0^\ell p[(f')^2-rf^2]\dd x=\int_0^\ell p w^2\left[\left(\frac fw\right)'\right]^2\dd x.
\end{equation}
If two boundary solutions existed, their difference $f$ would vanish at both endpoints and solve the homogeneous equation. Multiplying by $f$, integrating and using \eqref{eq:ground} forces $(f/w)'=0$; the endpoint condition then gives $f=0$.
\end{proof}
\begin{corollary}[Endpoint Wronskian and independent drift reconstruction]\label{cor:wronsk}
Under the preceding hypotheses, for any $\lambda\ge0$,
\begin{equation}\label{eq:wronsk}
 \frac{q_\lambda'(0)}{q_\lambda'(\ell)}=\exp\left(\int_0^\ell\frac{b(x)}{D(x)}\dd x\right).
\end{equation}
After recovering $D$, let $V_{\mathrm{phys}}=-\lambda_1/D-\frac12\mathcal S(q_{\lambda_1})$. The unique positive solution of
\begin{equation}\label{eq:phibvp}
 \phi''=V_{\mathrm{phys}}\phi,\qquad \phi(0)=1,\quad
 \phi(\ell)=\sqrt{q_{\lambda_1}'(0)/q_{\lambda_1}'(\ell)}
\end{equation}
gives the same drift by $b=2D\phi'/\phi$. In particular the endpoint-slope ratio is invariant in reset rate.
\end{corollary}
\begin{proof}
The Wronskian $W=u_-u_+'-u_-'u_+$ satisfies $W'=-gW$. At the endpoints $W(0)=-q_\lambda'(0)$ and $W(\ell)=-q_\lambda'(\ell)$, giving \eqref{eq:wronsk}. The gauge $\phi(x)=\exp(\frac12\int_0^x g)$ satisfies $\phi''=(g'/2+g^2/4)\phi=V_{\mathrm{phys}}\phi$ by the Schwarzian identity. The positive ground-state factorization gives uniqueness of the boundary solution; the drift formula is immediate.
\end{proof}
\begin{remark}[A rate-zero control]
When $\lambda_1=0$, $u_-+u_+=1$, so $w\equiv1$ and $b/D=-q_0''/q_0'$. If $D$ is already known, a single full positive-rate profile with \eqref{eq:wbvp} also reconstructs the drift without taking a third spatial derivative.
\end{remark}
\section{Conditioning, scaling, and numerical evidence}\label{sec:numerics}

\begin{lemma}[Scaling]\label{lem:scaling}
Let $q_\lambda[D,b;\ell]$ be the profile on $[0,\ell]$. For $\gamma>0$, $q_{\gamma\lambda}[\gamma D,\gamma b;\ell]=q_\lambda[D,b;\ell]$ (time change), and $q_\lambda[D,b;\ell](\ell z)=q_\lambda[D/\ell^2,b/\ell;1](z)$ (space change). Hence the conditioning of the two-rate problem depends on the rates only through $\lambda\ell^2/D$ for a given shape of $(D,b)$. 
\end{lemma}
\begin{proof}
$u_\pm$ solve $Du''+bu'=\lambda u$ with fixed boundary values; both changes of variable map this boundary problem to the stated one.
\end{proof}

\paragraph{Observations.}  Example $D=1+0.4\sin2\pi x+0.2\cos5x$, $b=0.8\cos3x-0.5$ on $[0,1]$, $\lambda_1=0$, $201$ observation points, independent Gaussian noise of standard deviation $\delta$ added to the values of both profiles (reproducible generator; see the README of the programs). Relative errors are root-mean-square on $[0.1,0.9]$, and $\bar D$ is the mean of $D$.
\begin{enumerate}[label=(\arabic*)]
\item \emph{Direct formula.} Smoothing spline (GCV) and the Schwarzian, $\lambda_2=20$: relative error of $D$ on $[0.1,0.9]$ is $4\cdot10^{-2}$ without noise, $0.3$ for $\delta=10^{-5}$, of order $10$ for $\delta\ge10^{-4}$.
\item \emph{Degree of ill-posedness.} Linearized response of the profiles to a mode $\sin(k\pi x)$ of $\log D$, after the best compensating change of $b$: decays like $k^{-3.2}$ to $k^{-3.7}$; for $b$, like $k^{-2.6}$ to $k^{-2.9}$ ($4\le k\le20$, $\lambda_2\le20$). This matches the formulas: $D$ needs three derivatives of $q$, $b$ two.
\item \emph{Finite-dimensional fit.} $\log D$ and $b$ in Chebyshev polynomials of degree $8$, nonlinear least squares (Levenberg--Marquardt from $\log D=0$, $b=0$) on both profiles, $\lambda_2=20$: error of $D$ $1.5\cdot10^{-2}$ without noise (truncation of the model), $1.7\cdot10^{-2}$ for $\delta=10^{-4}$ and $1.1\cdot10^{-1}$ for $\delta=10^{-3}$. About $100$ times more noise tolerated than with the direct formula.
\item \emph{Best second rate.} Smallest singular value of the Jacobian (8 modes of $\log D$, 8 of $b$): $\lambda_2=2,5,20,80,200,500,2000$ give $1.2\cdot10^{-4}$, $2.6\cdot10^{-4}$, $5.4\cdot10^{-4}$, $2.0\cdot10^{-4}$, $3.7\cdot10^{-5}$, $4.7\cdot10^{-6}$, $7.8\cdot10^{-7}$; for another $(D,b)$ the maximum is at $20$ to $40$. Small $\lambda_2$: the two profiles are too alike to separate $D$. Large $\lambda_2$: $D$ is well seen but $q_\lambda$ hardly depends on $b$; moreover the least-squares landscape develops local minima (for $\lambda_2=500$ the fit from $\log D=0$, $b=0$ does not reach the solution even without noise). By Lemma~\ref{lem:scaling}, the observed rule of thumb is $\lambda_2\ell^2/\bar D\approx20$, i.e.\ a typical distance $\sqrt{2\bar D/\lambda_2}\approx\ell/3$ travelled between resets; this is an observation for these examples, not a general prescription.
\end{enumerate}

\section{Interpretation and limitations}
The interval and endpoint-accumulation theorems give two conditions forcing reflection symmetry from restart neutrality. Neither conclusion follows from a single prescribed off-midpoint pair: the periodic-half construction works at every $0<y<L/2$. The first high-energy expansion retains the $o(k^{-1})$ remainder appropriate for an arbitrary $L^1$ potential. Pointwise continuity at the four specified locations controls the second high-energy term and yields the additional endpoint condition of Corollary~\ref{cor:necessary}; the unrestricted expansion to all orders remains formal. Exact two-rate reconstruction assumes continuum profiles and their derivatives; it does not establish stable inversion from finitely sampled noisy data. Spectral stability questions are separate from the exact identifiability statement~\cite{SavchukShkalikov}.

\section{Verification}\label{sec:verification}
The results are proved analytically under their stated hypotheses. Supplementary exact, interval-certified and numerical computations provide independent checks of selected identities, examples and reconstructions. The following table records the scope of these checks.

\begin{table}[htbp]\centering\scriptsize
\setlength{\tabcolsep}{3pt}
\renewcommand{\arraystretch}{1.18}
\begin{tabular}{>{\raggedright\arraybackslash}p{0.255\linewidth}>{\raggedright\arraybackslash}p{0.105\linewidth}>{\raggedright\arraybackslash}p{0.315\linewidth}>{\raggedright\arraybackslash}p{0.245\linewidth}}
\hline
Result & Proof & Exact / certified / numerical checks & Lean 4 \\
\hline
Renewal identities & Analytic & --- & Proposition~\ref{prop:pairform} \\
Reflection duality & Analytic & Exact P01--P02; certified C01 & Lemmas~\ref{lem:pos}, \ref{lem:rev}; Theorem~\ref{thm:real}(i),(iii)$^{\ast}$ \\
Sturm asymptotics & Analytic & Symbolic P03, P07; certified C02 & --- \\
Neutral-pair rigidity & Analytic & Exact P06 & --- \\
Non-symmetric neutral pairs & Analytic & Exact P05, P08, P09, P11, P12; certified P10; Python & Lemmas~\ref{lem:revmat}, \ref{lem:matrix}; Theorem~\ref{thm:families}(a),(c) and sufficiency in (b) \\
Midpoint Dirichlet / Weyl distinction & Analytic & Exact symbolic P14 & --- \\
Jacobi neutrality and recovery & Analytic & Exact rational P13; SymPy & Lemma~\ref{lem:jacobi}, Theorem~\ref{thm:disc}$^{\ast}$, Proposition~\ref{prop:three} \\
Two-rate reconstruction & Analytic & Exact P04; floating point C03; Python & Proposition~\ref{prop:schwarz}; \eqref{eq:wbvp}--\eqref{eq:brec} \\
Conditioning and rate selection & Numerical only & Floating point E01--E03 & --- \\
\hline
\end{tabular}
\caption{Scope of the mathematical proofs and supplementary checks. P, C and E refer to the programs in the folders \texttt{symbolic}, \texttt{checks} and \texttt{exploration} of the C++17 program suite. $^{\ast}$In the form stated below the table.}
\label{tab:verification}
\end{table}

The C++17 programs use the Herramientas library (version 0.14.0), including GiNaC for exact symbolic computation and FLINT/Arb for certified interval arithmetic. The collection comprises fourteen numbered programs, three checks and three numerical explorations; the accompanying README identifies the purpose and type of each calculation. A separate Python implementation using SymPy and mpmath checks the $\mathrm{SL}_2$ identities of Lemma~\ref{lem:matrix} and Theorem~\ref{thm:families}(b) on $1080$ exact rational instances. It tests the periodic-half family at five rates in 82-digit arithmetic and the Schwarzian and endpoint-Wronskian formulas of Section~\ref{sec:inverse}. The SymPy calculations also check the Cramer formulas and four-channel reconstruction on six rational chains ($a=4,\ldots,9$), and the two-piece Weyl identity of Example~\ref{ex:weyl}.

The Lean~4/Mathlib project \texttt{LeanPaperB1} formalizes the results recorded in the table, in the following form. In the continuous part, solutions are taken in the Carath\'eodory sense, as in Section~\ref{sec:framework}, with $V$ only integrable; uniqueness for the initial value problem and the existence of fundamental solutions enter as explicit hypotheses. This covers Lemmas~\ref{lem:pos} and~\ref{lem:rev}, and Theorem~\ref{thm:real}(i) for $V=\wt V$ almost everywhere; for (iii) it is the identity $q_\lambda(\mu)=C^\ast$ for every $\lambda>0$, whose reading for every renewal clock is Proposition~\ref{prop:laws}. For Lemma~\ref{lem:revmat}, the transfer matrices are shown to have unit determinant, to satisfy the product rule for concatenated pieces and the reversal formula; together with Lemma~\ref{lem:matrix}, this gives Theorem~\ref{thm:families} for the potential itself in its three cases, except the necessity in part~(b), which rests on the high-energy bound of Lemma~\ref{lem:volterra}. In the discrete part, the expansion of tridiagonal determinants, the invertibility of $I-sP$ and the formulas of Lemma~\ref{lem:jacobi} are proved, the channels entering through the first-step systems; Theorem~\ref{thm:disc} is proved in the form ``$U_+/U_-$ is constant on $(0,1)$ if and only if $a=2z$ and $d_L(z-1;s)=d_R(z+1;s)$'', whose reading for every renewal clock is Lemma~\ref{lem:discclock}, and Proposition~\ref{prop:three} follows for the chain itself. Proposition~\ref{prop:pairform} is formalized as stated, with the positivity of Lemma~\ref{lem:pos} as a hypothesis. In Section~\ref{sec:inverse}, the Schwarzian identity of Proposition~\ref{prop:schwarz} is proved for the ratio of two solutions of $y''+Py'+Qy=0$, its derivatives up to order three being computed, together with its consequences in that proposition; the derivation of \eqref{eq:wbvp} and \eqref{eq:brec} in the proof of Theorem~\ref{thm:two} is proved pointwise. Everything else is not formalized: the renewal arguments (Propositions~\ref{prop:single} and~\ref{prop:laws}, Lemma~\ref{lem:discclock}) and the first-step identification of the channels; Lemmas~\ref{lem:volterra}--\ref{lem:weylunique}; Section~\ref{sec:asym}; Theorem~\ref{thm:real}(ii), Theorems~\ref{thm:weights}--\ref{thm:endpoint}, the realization by a drift at the end of Theorem~\ref{thm:families} and Corollary~\ref{cor:necessary}; Section~\ref{subsec:midweyl}; the existence, positivity and uniqueness argument of Theorem~\ref{thm:two} and Corollary~\ref{cor:wronsk}; and Section~\ref{sec:numerics}. The Lean proofs use only the standard axioms. No computation replaces an analytic proof. The C++ programs and the Lean project are archived at Zenodo, \href{https://doi.org/10.5281/zenodo.24402657}{doi:10.5281/zenodo.24402657}.

\section*{Acknowledgements}
Both authors used AI tools in this work: T.~Newton used OpenAI ChatGPT and J.~A. Vega Coso used Claude (Anthropic). The tools aided exploration, symbolic and formal verification (including code for exact computations and for the Lean formalization), numerical exploration, and the preparation of the text. The algorithms, research direction and mathematical content remain the authors' responsibility. The theorems were accepted only on the basis of proofs, supported by exact, interval-certified or formal checks; the numerical observations of Section~\ref{sec:numerics} are reproducible floating-point computations and are not claimed as theorems.

\end{document}